\iffalse
\documentclass[12pt,a4paper]{article}
\usepackage{setspace}
\let\appendices=\appendix
\else
\documentclass[conference]{IEEEtran}
\fi
  \usepackage{cite}
\usepackage{algorithm}

\usepackage{algpseudocode}

  \usepackage[pdftex]{graphicx}
\usepackage{amsmath,amssymb}
\usepackage{bbm}
\usepackage{mathtools}
\usepackage{hyperref}

\def\etal{\emph{et al.}}

\begin{document}

\title{On Stochastic Rounding with Few Random Bits}

\author{
\IEEEauthorblockN{Andrew Fitzgibbon}
\IEEEauthorblockA{Graphcore\\Cambridge, UK\\Email: awf@ieee.org}
\and
\IEEEauthorblockN{Stephen Felix}
\IEEEauthorblockA{Graphcore\\Bristol, UK\\Email: Stephen.Felix@graphcore.ai}
}

\date{}


%


\maketitle

\begin{abstract}
Large-scale numerical computations make increasing use of low-precision (LP) floating point formats and mixed precision arithmetic, which can be enhanced by the technique of \emph{stochastic rounding} (SR), that is, rounding an intermediate high-precision value up or down randomly as a function of the value's distance to the two rounding candidates.  Stochastic rounding requires, in addition to the high-precision input value, a source of random bits.  As the provision of high-quality random bits is an additional computational cost, it is of interest to require as few bits as possible while maintaining the desirable properties of SR in a given computation, or computational domain.   This paper examines a number of possible implementations of few-bit stochastic rounding (FBSR), and shows how several natural implementations can introduce sometimes significant bias into the rounding process, which are not present in the case of infinite-bit, infinite-precision examinations of these implementations.  The paper explores the impact of these biases in machine learning examples, and hence opens another class of configuration parameters of which practitioners should be aware when developing or adopting low-precision floating point.
Code is available at \url{http://github.com/graphcore-research/arith25-stochastic-rounding}.
\end{abstract}


%

\section{Introduction}
Large-scale numerical computations make increasing use of low-precision (LP) floating point formats, both for storage and computation.
As compared to 32 and 64 bit formats, narrow storage formats (16, 8, and even fewer bits per element) allow larger arrays to be stored in high-bandwidth memory, close to accelerated arithmetic units; and allow more elements to be transported per second across communications links, whether within or between compute nodes.  Narrow-format arithmetic can be implemented more efficiently, both in hardware and software, increasing the speed of computation.  Of course, narrow formats imply decreased accuracy, requiring computations to be re-engineered in order to compensate.  In areas such as machine learning, the training of large models has demonstrated that this re-engineering can yield significant improvement in the scale and speed at which a computational goal can be accomplished, while yielding models of equivalent or greater accuracy at reduced cost.  In some instances, the low-precision models are superior in the sense that they can be deployed on lower-power devices after training.

This re-engineering, however, can prove difficult or impossible at lower format widths.  There is, after all, a considerable loss of precision and dynamic range in LP formats, giving rise to the following challenges: careful dynamic range management, judicious deployment of mixed low and standard precision, a re-exploration of ``hyper-parameter'' settings (such as step sizes, number of iterations, etc).
These challenges can be mitigated by the use of stochastic rounding (SR), as has been shown in many application domains~\cite{essam17,noune2022,paxton2022climate,zhang2022fast,croci2023effects}.

\begin{figure*}
  \centering
  \includegraphics[width=\linewidth]{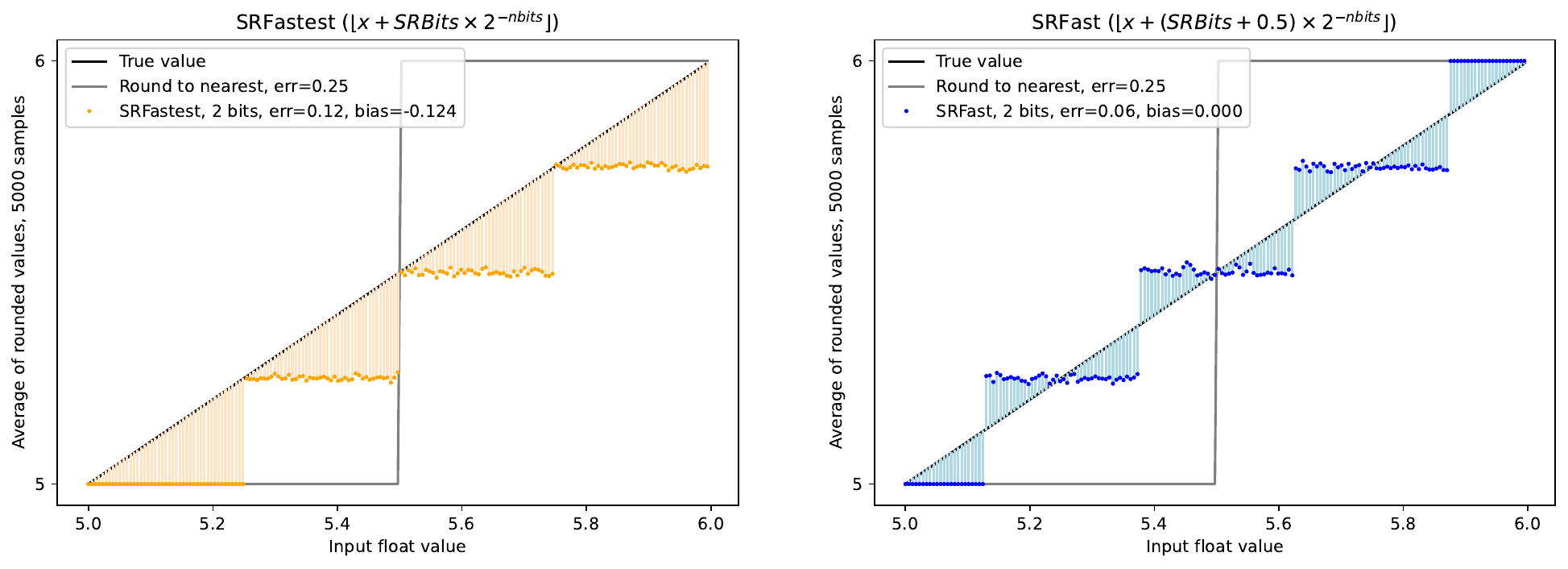}
  \caption{Bias in stochastic rounding with finite randomness.
  For each $X$ value, 5000 samples are rounded with 2 bits of
  randomness,
  and the mean of the rounded values is plotted.
  (Left, SRFF): Mean rounded values are below the line $Y=X$, indicating bias.
  The computed bias of $-0.124 \approx-2^{-3}$, as predicted by the calculation of \S\ref{sec:bias-srff}.
  (Right, SRF): Mean rounded values are symmetric around $Y=X$, with bias $\approx$ 0,
  following \S\ref{sec:unbiased-srf}.
  }
  \label{fig:bias1}
\end{figure*}

The focus of this paper is on implementations of SR which use small numbers of random bits, which we call \emph{few-bit stochastic rounding} (FBSR).  In particular, we consider the case where the number of random bits used is smaller than the difference in precision between the values to be rounded and the target precision.  We illustrate (see Figure~\ref{fig:bias1}) that natural implementations may have sometimes significant biases, confounding any experimental investigations into FBSR.  Mitigations for these biases are presented, and empirical results are presented to argue for the effectiveness of these mitigations, and hence for the utility of FBSR.

\subsection{Related work}
Croci \etal~\cite{croci2022stochastic} provide an excellent survey of the theory and applications of SR in a variety of fields including machine learning and quantum computing.  In almost all current implementations, whether hardware or software, the number of random bits used equals the number of trailing significand bits, so are not in the few-bit regime described in this paper.  For future implementations, however, it is expected that our findings will be useful in making design decisions.

This paper's primary contribution is a description of, and remediation of, the bias that some implementations of FBSR may introduce.  This bias has been analysed theoretically in recent work by El Arar \etal~\cite{elarar23,elarar24}, but no mitigation is proposed.  The current paper proposes two alternative mitigations, and provides a more application-focused derivation of the bias computation which may prove useful to some readers.  Xia \etal~\cite{xia20} also consider bias in SR, but with the goal of introducing bias in order to reduce variance.   Given current usage of SR in large-scale machine learning applications where millions or billions of roundings occur in one training iteration, it is suggested that bias is of greater concern than variance in current applications.

The paper does not discuss the quality of the supplied bits, merely their quantity.
Some recent work~\cite{yuan2022you} suggests that lower quality sources (e.g.\ nearby truncations) may serve effectively in practice.  It is thus conceivable that more bits from a lower-quality source might allow a more efficient implementation than a bias-free implementation with fewer high-quality bits, although the extent of the bias is such that this appears unlikely, and the conflations involved in using lower-quality sources might render theoretical analysis difficult or impossible.

As noted by several authors~\cite{essam17,noune2022}, low-precision formats generally make use of a \emph{scaled tensor} implementation, where blocks of low-precision values are accompanied by a scale factor, which may be a power of two
(i.e.\ a pure-exponent format, or another floating point format).
This somewhat expands the applicability of the approach described in this paper, as even if intermediate values are computed in mid-precision formats such as BFloat16, the scaled values may have larger effective precision, making the corrections in this paper more important.

\subsection{Contributions and limitations}
\noindent
For the reader's convenience, the contributions and limitations of this work are outlined:
\begin{itemize}
\item A simple but general mathematical framework for the description and analysis of FBSR, with a particular focus on practical implementation;
\item A definition of several FBSR methods within this framework;
\item Analysis of the bias properties of these methods, revealing biases that are not, to our knowledge, described in the existing literature;
\item Experimental investigations of the practical impact of the discovered biases on representative computations from machine learning.
\end{itemize}

\noindent
Primary limitations are:
\begin{itemize}
\item
The experimental investigations are at a scale which is sufficient to confirm that the biases can have an effect on real-world use cases, but these conclusions may not extend to larger-scale deployments.
\item
The costs, in terms of area and power, of any hardware implementation are not discussed.  This is at least in part because so many area/power decision tradeoffs exist even when the parameters are fixed.  It is hoped that the paper's investigations are of utility in making such tradeoffs in any future hardware designs.
\end{itemize}

\def\ebase{{e_0}}


\begin{figure*}
  \centering
  \includegraphics[width=\linewidth]{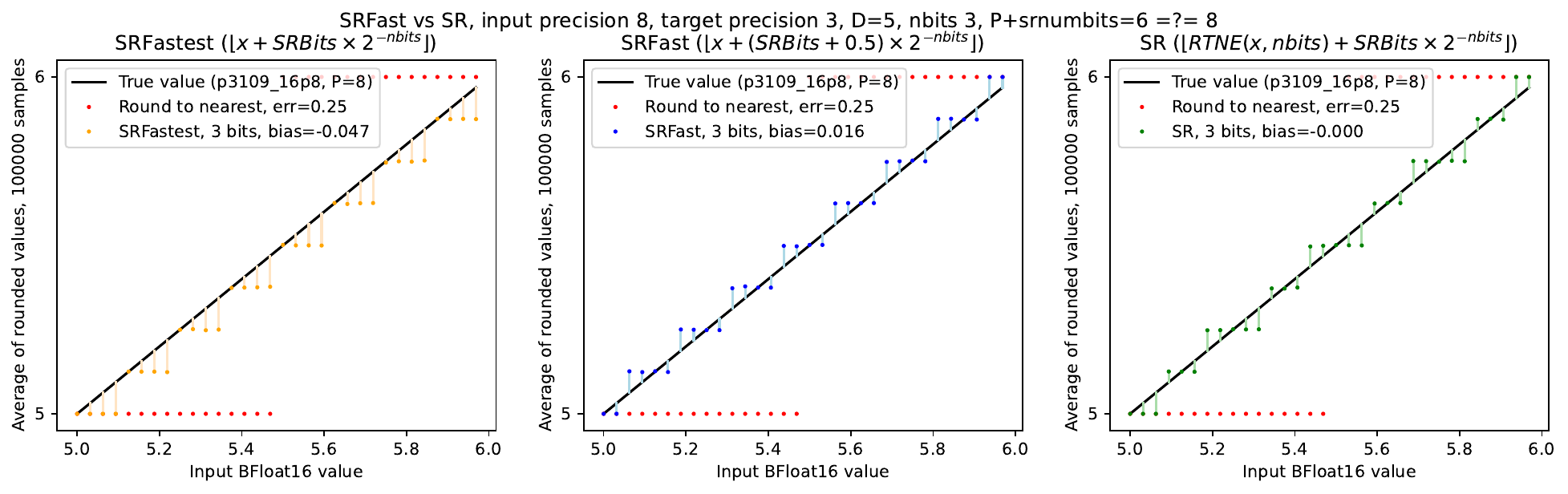}
  \caption{Bias in stochastic rounding with finite randomness, finite-precision inputs.
  For each $X$ value in BFloat16, 100,000 samples are rounded with 2 bits of
  randomness,
  and the mean of the rounded values is plotted.
  (Left, SRFF): Mean rounded values are below the line $Y=X$, indicating bias.
  The computed bias is $\approx$ -0.046875, as predicted by the calculation of \S\ref{sec:fbbias-srff} and \S\ref{asec:fbbias-srff}.
  (Middle, SRF): Mean rounded values are not symmetric around $Y=X$, with bias $\approx$ 0.015625, following \S\ref{asec:fbbias-srf}.
  (Right, SRC): Mean rounded values are symmetric around $Y=X$, with empirical bias $\approx$ 0.000, following the correction of \S\ref{sec:fbbias-srff}.
  }
  \label{fig:bias2}
\end{figure*}

\def\comment#1{\text{#1}}
\section{Background and notation}
We consider binary floating point numbers represented by the triple $(s,e,t) \in \mathbb N^3$, representing sign, (biased) exponent and the fractional part (trailing bits) of the significand.
A floating point \emph{format} is parameterized by width in bits $K$, precision in bits $P$, and exponent bias $B$.  We shall also make use of~$M$, the largest finite value in the format's value set (from which it assumed symmetrically that $-M$ is the smallest finite value).
The real value encoded by $(s,e,t)$ in format $(K,P,B;M)$ is given by the function
$$
F(s,e,t; P,B)
$$
defined as
\begin{equation}
(-1)^s \times
    \begin{cases}
        (0 + t \times 2^{1-P}) \times 2^{1 - B} & \text{if~} e = 0\\
        (1 + t \times 2^{1-P}) \times 2^{e - B} & \text{if~} e > 0
   \end{cases}
\end{equation}
where the $e=0$ case encodes subnormal numbers.  In this definition, the precision~$P$ is one more than the number of mantissa bits that need to be encoded in a packed binary representation of the number.
Formats may additionally encode infinities, negative zero, and not-a-number (NaN) values, which shall largely not feature in the discussions of this paper (in implementations, it is assumed that values outside the representable finite float range are rounded deterministically as in IEEE-754 round to nearest).

\def\floor#1{\lfloor #1 \rfloor}
\def\R#1{\tilde{#1}}
\subsection{Rounding}
A strictly positive real value $\R x > 0$ is rounded as follows.
We first compute the real-valued significand $\R \sigma \in [2^{P-1}, 2^P)$, which is then rounded to integer $\sigma \in \{2^{P-1},...,2^P\}$, of which the trailing bits are encoded (incrementing the exponent if $\sigma=2^P$).
The initial exponent is computed as follows, allowing for subnormals:
\begin{equation}
 \R e = \max(\floor{\log_2(\R x)}, 1-B)
\end{equation}
from which the (scaled) significand is
\begin{equation}
  \R \sigma = \R x \times 2^{-e} \times 2^{P-1}
\end{equation}
Then, rounding to integer with ties away from zero
is represented by
\def\aIf{\mathrel{\mathsf{if}}}
\def\aOr{\mathrel{\mathsf{or}}}
\def\aAnd{\mathrel{\mathsf{and}}}
\def\IsEven{\mathop{\mathsf{IsEven}}}
\def\rsig{\R\sigma}
\def\rdel{\R\delta}
\def\floors{\floor{\rsig}}
\def\op#1{\operatorname{\mathsf{#1}}}
\def\roundawayonly{\op{R}}
\def\roundaway#1{\roundawayonly_{\text{#1}}}
\def\indicator#1{\mathbbm{1}\!\left[{#1}\right]}
\begin{eqnarray*}
 \rdel & = & \rsig - \floors\\
 \sigma & = &
 \begin{cases}
     \floors & \aIf \rdel < 0.5\\
     \floors+1 & \aIf \rdel \ge 0.5
 \end{cases}
\end{eqnarray*}
which may be written compactly in terms of a predicate
\begin{eqnarray}
\roundaway{TA}(\rdel) = \rdel \ge 0.5
\label{eq:predTA}
\end{eqnarray}
whence we can replace the above with
\begin{eqnarray}
 \sigma = \floors + \indicator{\roundaway{RoundMode}(\rsig - \floors)}
 \label{eq:pred}
\end{eqnarray}
For ties to nearest even or odd (TNE, TNO), we may define an alternative $\roundaway{RoundMode}$,
which has access to $\rsig$ in order to determine the direction in which to resolve ties.
For example,
$$
\roundaway{TNE}(\rsig, \rdel)
=
\rdel > 0.5 \aOr \bigl(\rdel = 0.5 \aAnd \IsEven(\floors+1)\bigr)
$$
where $\IsEven(i)$ is true iff integer $i$ is even.
In practice, implementations may use comparison as indicated above, or an alternative implementation using addition, analogous to $\rdel + 0.5 \ge 1.0$.  The conclusions of this paper apply to both implementations, with possible inversion of the direction of any bias.

\subsection{Stochastic Rounding}
With the above notation, stochastic rounding is readily defined
\[
\roundaway{SR}(\rdel) = \rdel + \R n \ge 1
\]
where the ``noise'' $\R n$ is a random variable drawn from a uniform distribution between~0 and~1.
As mentioned in the introduction, the quality of the random number generator is orthogonal to the concerns of this paper.
In practical implementations of SR, the random numbers will be supplied as a finite sequence of random bits, denoted by an integer $n$ such that
\[
0 \le n<2^N,
\]
where $N$ is the number of bits in $n$.
Random bits are passed to the $\roundawayonly$ function, so that a natural analogue of the infinite precision definition is given by
\[
\roundaway{SRFF}(\rdel, n) = \rdel + n \times 2^{-N} \ge 1
\]
As will be explored in the following, this definition proves to be biased, and variations which mitigate this bias shall be defined.  One which is of common practical utility is the form we shall name ``SRF''\footnote{The naming scheme encodes speed (and concomitant reduced accuracy): corrected (SRC), fast (SERF), faster (SRFF).}:
\[
\roundaway{SRF}(\rdel, n) = \rdel + \left(n + \frac12\right) \times 2^{-N} \ge 1
\]

\section{Bias}
\label{sec:bias-calcs}
\def\pfam{\mathcal{P}}
A property that one might require of any rounding scheme is that it is \emph{unbiased}, or that the expected difference between inputs and outputs is zero.
Written loosely, this might be
\[
\mathbb E[x-\op{round}(x)] = 0.
\]
As written, this definition is somewhat vacuous: the expectation can only be over values of $x$, and hence depends on a choice of prior probability distribution $p(x)$.
Writing more explicitly, the requirement is
\[
\mathbb E_{x\sim p(x)}[x-\op{round}(x)] = 0.
\]
This prior should not, for example, simply be the set of real values which map to the format's finite range, e.g.\ $p(x)$ being the uniform distribution on $[-M,M]$.  Under such a distribution, always rounding toward zero is unbiased, as the negative values will cancel the positives.

More meaningfully, the expectation may be over a family of {\em test distributions} $\pfam$, and we require unbiased rounding for all members of the family:
\[
\kern -3em \forall p \in \pfam: \quad \mathbb E_{x\sim p(x)}[x-\op{round}(x)] = 0.
\]
In this paper, we choose to require that values between each pair of floats are unbiased on average, i.e., if $\mathcal U(l,h)$ is the uniform distribution between $l$ and $h$, and $\mathcal F$ is the set of finite floats in a format, then a family of test distributions is
\begin{equation}
  \pfam = \bigl\{\mathcal U(f, \operatorname{succ}(f)), f \in \mathcal F - \{M\}\bigr\}
  \label{eq:pfam}
\end{equation}
hence we require that
\[
\forall f\in \mathcal F, f < M: \quad \int_f^{succ(f)} (x-\op{round}(x)) \mathrm dx = 0
\]
In terms of a rounding predicate $\roundawayonly$, as in (\ref{eq:pred}), unbiased rounding then means
\begin{eqnarray*}
  \int_0^1 \indicator{\roundawayonly(\rdel)} - \rdel \mathrm d\rdel & = & 0 \\
\end{eqnarray*}
or, equivalently
\begin{eqnarray*}
  \int_0^1 \indicator{\roundawayonly(\rdel)} d\rdel - \frac12 & = & 0
\end{eqnarray*}

\subsection{Stochastic rounding}
For SR, we must additionally compute the expectation over the random variable $\R n$:
\[
\kern -3em
\forall p \in \pfam: \quad
\mathbb E_{\R n}\bigl[ \mathbb E_{x\sim p(x)}[x-\op{round}(x)]\bigr] = 0.
\]
In the infinite-precision case, the family $\pfam$ can be delta functions at all $x \in [-M,M]$, with zero bias at all points
\[
  \forall x \in [-M,M]: \quad
  \int_0^1 [x- \op{round}(x,\R n)] \mathrm d\R n = 0
\]
Equivalently, this is the requirement that
\[
  \forall \R \delta \in [0,1]: \quad
  \int_0^1 \indicator{\R\delta + \R n \ge 1}  \mathrm d\R n = \R \delta
\]
which is straightforwardly true as for all $\R \delta \in [0,1]$ we have
\begin{eqnarray*}
  \int_0^1 \indicator{\R\delta + \R n \ge 1}  \mathrm d\R n &=& \\
  \int_0^1 \indicator{\R n \ge 1 - \R\delta}  \mathrm d\R n &=& \\
  \int_{1 - \R\delta}^1 1 \mathrm d\R n &=& \\
  1 - (1 - \R\delta) &=& \R \delta
\end{eqnarray*}

\subsection{Stochastic rounding, finite randomness}
In the finite-randomness case, the test distributions must revert to the family of~(\ref{eq:pfam}), i.e. uniform between each pair of floats, so the unbiasedness condition is
\[
\mathbb E_n\left[\int_0^1 \indicator{\roundawayonly(x,n)} \mathrm dx\right]  - \frac12 = 0
\]
which is the discrete sum
\[
2^{-N}\sum_{n=0}^{2^N-1} \int_0^1 \indicator{\roundawayonly(x,n)} \mathrm dx  - \frac12 = 0
\]

\subsection{Bias of SRFF}
\label{sec:bias-srff}
Let us now compute the bias of the SRFF mode, for a given number of bits $N$.
Appendix~\ref{app:srff} contains the complete calculation, key steps of which are as follows:
\begin{eqnarray*}
2^{-N}\sum_{n=0}^{2^N-1} \left(\int_0^1 \indicator{\roundaway{SRFF}(x,n)} \mathrm dx  - \frac12 \right)&=&\\
2^{-N}\sum_{n=0}^{2^N-1} \int_0^1
\indicator{x + n\times2^{-N} \ge 1} \mathrm dx- \frac12 &=&\\
2^{-N}\sum_{n=0}^{2^N-1} \int_{1-n\times2^{-N}}^1
 \mathrm dx - \frac12&=& \\
-2^{-(N+1)}&\ne& 0
\end{eqnarray*}
Hence SRFF is biased, with bias decreasing with increased number of random bits.  As shown later in experiments, even with 4 bits of randomness, the bias can significantly affect deep learning training.

\subsection{Unbiasedness of SRF}
\label{sec:unbiased-srf}
A similar calculation for SRF mode becomes (again see appendix for detail)
\begin{eqnarray*}
2^{-N}\sum_{n=0}^{2^N-1} \int_0^1 \indicator{\roundaway{SRF}(x,n)} \mathrm dx  - \frac12 &=&\\
2^{-N}\sum_{n=0}^{2^N-1} \int_0^1
\indicator{x > \left(n+\frac12\right)\times2^{-N}} \mathrm dx- \frac12 &=& 0
\end{eqnarray*}
showing that SRF is unbiased, in that the expectation computed over all real values between two floats is zero.

\subsection{Bias, finite precision inputs}
\label{sec:fbbias-srff}
The above calculations assume that all real numbers are equally likely inputs to rounding.
In real-world systems, the inputs will be available at some precision $Q>P$, which we shall write in terms of the number of excess bits $D$, i.e.~$Q = P + D$.  When $D$ is large, for example, when converting from 32-bit or higher formats, we might expect the calculations above for infinite-precision inputs to pertain.
However it is not unusual for the source format to be relatively close in precision to the target.  For example, BFloat16 has a precision of~8 and might be rounded to an E5M2 format with a precision of~3, in which case $D=5$.

In such cases, the integral also becomes a discrete sum over the incoming values between~0 and~1.  For simplicity, we consider initially only the case where the incoming range does not wrap binades in our target interval.  Hence there are $2^D$ values (excluding 1).  Then the bias of SRFF becomes
\begin{multline}
  bias_{SRFF,D} =\\
  2^{-D} \sum_{i=0}^{2^D-1} \left( 2^{-N}\sum_{n=0}^{2^N-1}
    \indicator{\roundaway{SRFF}(2^{-D} i,n)} - 2^{-D} i
    \right)
\end{multline}
As derived (with more verbose annotations) in the appendix, this calculation yields the bound (tight for $N \le D$)
\begin{eqnarray*}
  bias_{SRFF,D} \le\frac{2^{-D} - 2^{-N}}2  \label{eq:fbsrff}
\end{eqnarray*}
from which it is clear that the bias is zero for $N=D$.

We note that $N=D$ is the case commonly implemented in today's hardware, 
however, as noted in the introduction, 
the research literature is increasingly exploring $N < D$, 
in which case correction of the bias will be required.

It is not sufficient to use SRF, which also has a bias, as shown in \S\ref{asec:fbbias-srf} bounded by (and again, with the bound being tight for $N<=D$)
\begin{eqnarray*}
  bias_{SRF,D} \le 2^{-(D+1)}
\end{eqnarray*}

\subsection{Bias-corrected few-bit SR (SRC)}
There is, however, a simple correction: first deterministically round the incoming precision-$Q$ values to precision $Q' = P + N$ then apply SRFF.
This reduces the problem to the $N=D'$ case, which is unbiased for SRFF as described above.  More formally:
\begin{eqnarray}
  \roundaway{SRC}(\R\delta,n) = \roundaway{SRFF}(\text{Round}(\R\delta \times 2^N) \times 2^{-N}, n)
\end{eqnarray}
where $\text{Round}$ is any unbiased rounding scheme, e.g.\ round to nearest, with ties to even or odd. 
Despite its simplicity, we are not aware of this correction having been described in any existing discussion of FBSR.  As shall be shown in the sequel, empirical results support the theoretical calculations above.

\begin{algorithm}[t]
  \begin{algorithmic}
    \State $w_8 : F8 \gets$ initializers
    \While{not done}
      \State $w_{16} : F16 \gets w_8$
      \State $g_{16} : F16 \gets \nabla \text{loss}(w_{16})$
      \State $u_{16} : F16 \gets $ AdamW update in F16
      \State $w_{8} \gets \text{RoundTo8}(w_{16} + u_{16})$ \Comment{SR applied here}
    \EndWhile
    \end{algorithmic}
  \caption{Quantization-aware training of an 8-bit (F8) model using high-precision (F16) gradients.}
  \label{alg:qat}
  \end{algorithm}
  
\begin{figure}
  \centering
  \includegraphics[width=\linewidth]{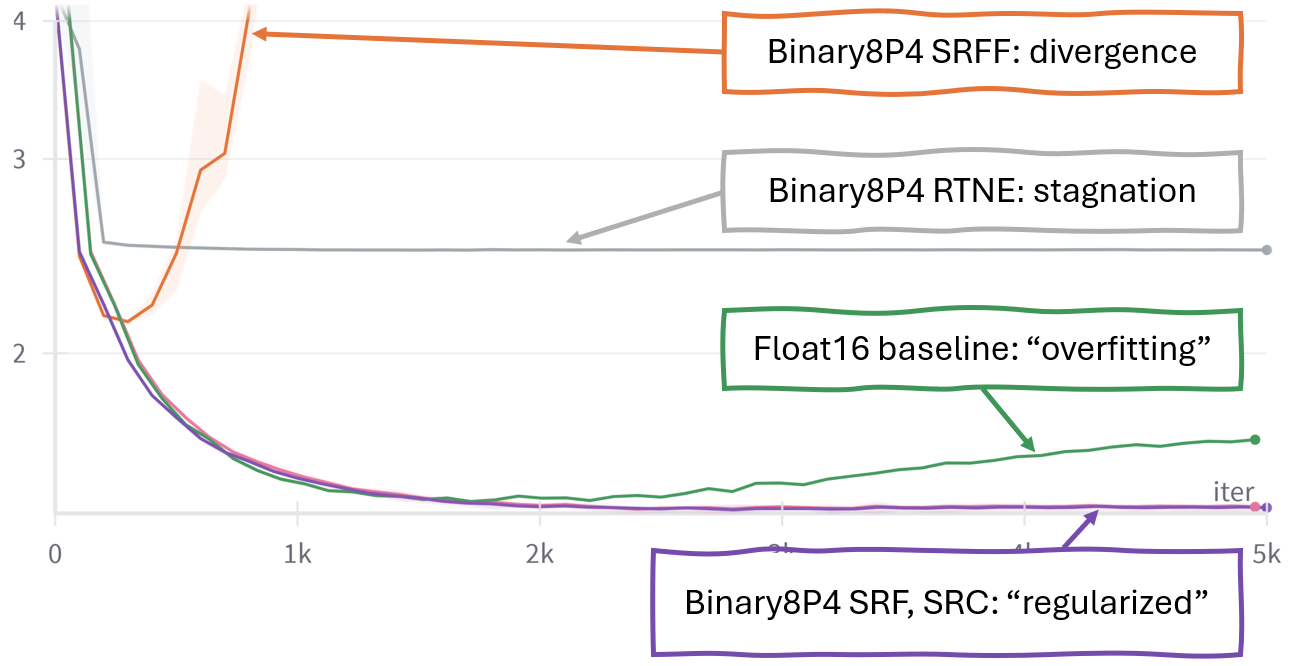}
  \caption{Experiments on language model training (nanoGPT, small model).
  Binary16 baseline: achieves minimum validation loss at about 1900 iterations, after which training loss (not shown) reduces, but validation loss increases.
  Binary8P4 RTNE: initial loss reduction followed by ``stagnation'';
  SRFF: initial loss reduction followed by ``divergence'';
  SRF, SRC: convergence to stable values.
  Although SRF and SRC converge to better validation errors than Binary16, this regularization behaviour is expected to apply only to such small-scale experiments, and should not be taken as an indication of any superiority of Binary8+SR over Binary16.
  }
  \label{fig:expt}
\end{figure}

\begin{figure*}[t]
  \centering
  \includegraphics[width=\linewidth]{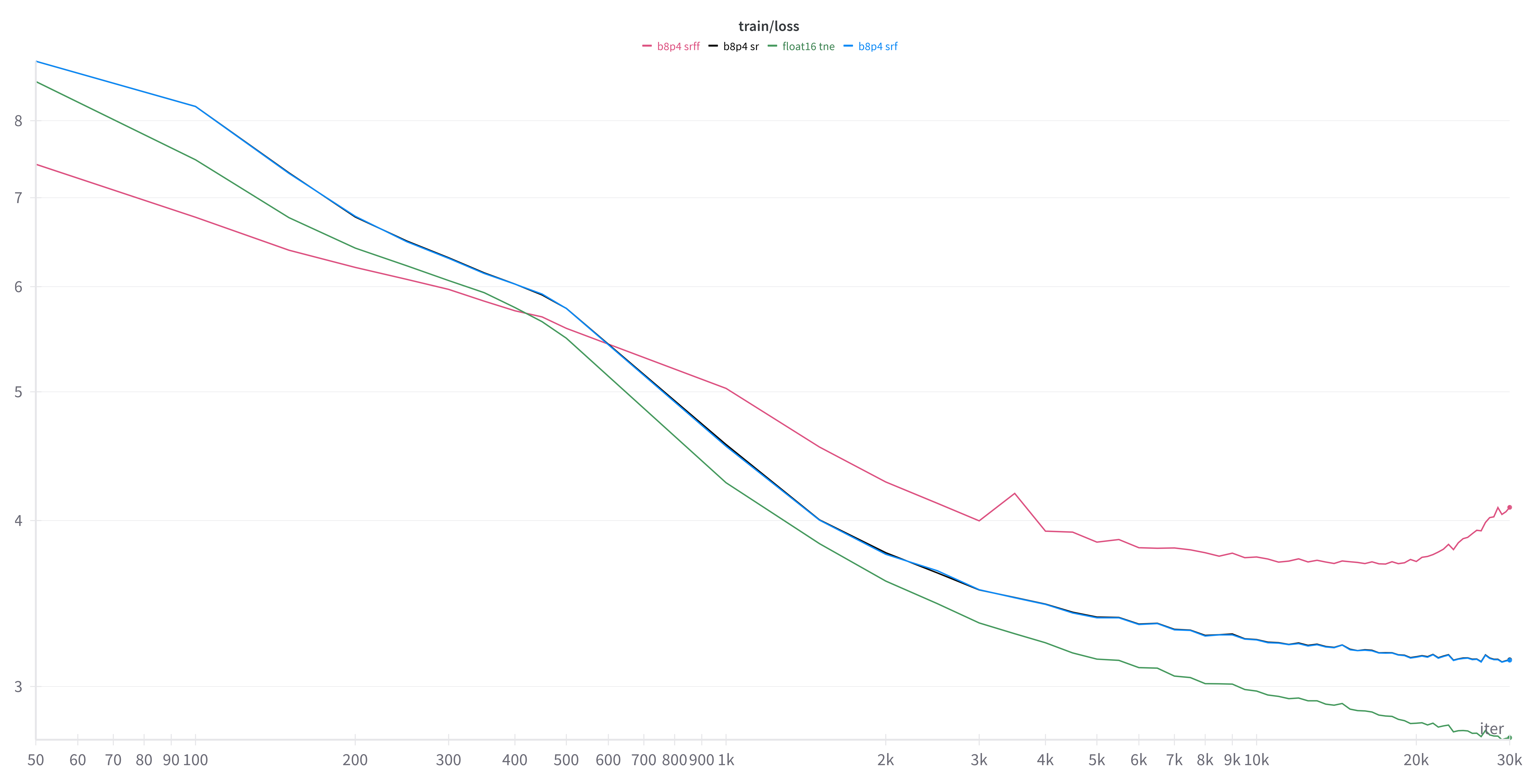}
  \caption{Experiments on language model training (nanoGPT, GPT-2 (350M params)).
  Binary16 baseline: validation loss follows training loss (not shown) very closely.
  Binary8P4 SRFF: initial loss reduction followed by ``divergence'';
  SRF, SRC: convergence to stable values.
  In this case SRF and SRC converge to better validation errors than SRFF,
  but not to the same loss as Binary16.
  }
  \label{fig:gpt2}
\end{figure*}

  \section{Experiments}
A current interest in the field of generative language modelling is in the training of models with low-precision weights, such that inference can be performed on low-power devices.
While low-precision weights can be obtained by post-training quantization (PTQ), that is, simply casting learned weights to the low-precision format, it is increasingly understood that performing quantization-aware training (QAT) can yield more accurate models~\cite{or24}.
To this end, we have performed a simple simulation of QAT whereby weights are rounded to the target precision at initialization and after every gradient update (or equivalently, before every gradient calculation).
A more precise definition in pseudocode is given at Algorithm~\ref{alg:qat}.

Algorithms were implemented on the {\tt nanoGPT} codebase~\cite{Karpathy2022}, with stochastic rounding implemented via the {\tt gfloat} library~\cite{Fitzgibbon24}.
The implementation did not make use of SR implemented in hardware, in order to be able to freely vary the SR implementation.
Code is available at \url{http://github.com/graphcore-research/arith25-stochastic-rounding}.

The default model was trained on the Shakespeare 
dataset.  At this small scale, the model is known to overfit training data, and hence validation loss is the figure of merit followed.
Training baselines were run with gradient updates computed in BFloat16, and Binary16, yielding essentially identical validation curves.
The target format was set to IEEE WG P3109 Binary8P4~\cite{p3109} format whose  precision corresponds to that of OCP E4M3 format.
For these experiments, the number of SR bits used was 3.  In principle, as gradient updates are in BFloat16, with a precision of~8, and the target has a precision of~4, the ideal number of SR bits to use is 4.

As illustrated in Figure~\ref{fig:expt}, rounding to nearest with ties to even (RTNE) succeeds in initially reducing the validation loss, but then ``stagnates'', where further gradient updates make no change to the weights and hence to the loss.
Rounding to Binary8P4 under the (biased) SRFF implementation initially decreases the loss, but soon ``diverges'', with the loss rapidly increasing.
Both bias-corrected implementations (SRF, SRC) decrease the loss, initially to the minimum value achieved by the Binary16 baseline, and then continue to slightly reduce loss.
However, as mentioned in the figure caption, this regularization behaviour is expected to apply only to such small-scale experiments, and should not be taken as an indication of any superiority of Binary8+SR over Binary16.
The model was then configured to GPT-2 settings (26 layers, 16 heads, embedding dimension 1024), amounting to 354 million parameters.  In these experiments, training was stopped at 30,000 iterations, meaning that we expect under-training, so validation loss follows training loss (this was observed).  The model was trained in Binary8P4 as above, under SRFF, SRF, SRC modes.
In this instance, SRFF attained a final loss value of 4.06, compared to 3.14 for each of SRC, SRF, and 2.74 for Binary16.  It should be noted that the training dynamics (and hence the optimal hyperparameter settings) of this simple QAT method are almost certainly quite different from the Binary16 defaults, so it is not necessarily expected that the loss achieved by SR will match Binary16.  The key outcome of this experiment is the indication that the correction of SRFF's bias can avoid divergence.  Divergence was not as dramatic as for the smaller model, but was nevertheless significant, suggesting that alternatives to SRFF may be of value in future implementations.

\section{Conclusion}
This paper has illustrated the potential for systematic bias in some plausible implementations of few-bit stochastic rounding.
The mitigations of these biases are also presented, with different cost/accuracy tradeoffs.  Some experimental evidence is provided to suggest that these biases may prove significant in practice, and that further exploration into the costs of hardware or software implementations may be of value.
Source code implementing all experiments, and including symbolic algebra derivations of the various computations, will be made available.




\appendices

\section{Bias computations}
This section records a detailed derivation of the bias computations in \S\ref{sec:bias-calcs}.

\subsection{Bias of SRFF, infinite-precision inputs}
The calculation of \S\ref{sec:bias-srff} is expanded here:
\label{app:srff}
\begin{eqnarray*}
 &&
2^{-N}\sum_{n=0}^{2^N-1} \int_0^1 \indicator{\roundaway{SRFF}(x,n)} \mathrm dx  - \frac12\\
&=&
2^{-N}\sum_{n=0}^{2^N-1} \int_0^1
\indicator{x + n\times2^{-N} \ge 1} \mathrm dx- \frac12\\
&=&
2^{-N}\sum_{n=0}^{2^N-1} \int_{1-n\times2^{-N}}^1  \mathrm dx - \frac12\\
&=&
2^{-N}\sum_{n=0}^{2^N-1} n\times2^{-N} - \frac12 \\
&=&
2^{-2N}\sum_{n=0}^{2^N-1} n - \frac12 \\
&=&
2^{-2N}\frac{(2^N-1)(2^N)}2 - \frac12 \\
&=&
2^{-2N}\frac{2^{2N}-2^N}2 - \frac12  \\
&=&
\frac{1-2^{-N}}2 - \frac12  \\
&=&
-\frac{2^{-N}}2  \\
&=&
-2^{-(N+1)} \ne 0
\end{eqnarray*}

\subsection{Bias of SRF, infinite-precision inputs}
\label{app:srf}
The calculation in \S\ref{sec:unbiased-srf} is expanded as follows:
\begin{equation*}
bias_{SRF} = I - \frac12  
\end{equation*}
where
\begin{eqnarray*}
I &=&
2^{-N}\sum_{n=0}^{2^N-1} \int_0^1 \indicator{\roundaway{SRF}(x,n)} \mathrm dx  \\
&=&
2^{-N}\sum_{n=0}^{2^N-1} \int_0^1
\indicator{x + \left(n+\frac12\right)\times2^{-N} \ge 1} \mathrm dx\\
&=&
2^{-N}\sum_{n=0}^{2^N-1} \int_{1-\left(n+\frac12\right)\times2^{-N}}^1  \mathrm dx  \\
&=&
2^{-N}\sum_{n=0}^{2^N-1} \left(\left(n+\frac12\right)\times2^{-N}\right)
\end{eqnarray*}
so
\begin{eqnarray*}
  bias_{SRF} &=& 2^{-N}\sum_{n=0}^{2^N-1} \left(n\times2^{-N} + \frac12\times2^{-N}\right) -\frac12\\
&=&
\text{bias}_{\text{SRFF}} + 2^{-N}\sum_{n=0}^{2^N-1} \left(\frac12\times2^{-N}\right)\\
&=&
-2^{-(N+1)} + 2^{-(N+1)} = 0
\end{eqnarray*}
\subsection{Bias of SRFF, finite-precision inputs}
\label{asec:fbbias-srff}
\def\ceil#1{\left\lceil#1\right\rceil}
\def\floor#1{\left\lfloor#1\right\rfloor}
The bias at a single $x$ value is
\begin{eqnarray*}
  bias_{SRFF}(x) \kern -10mm && \\
  &=&
  2^{-N}\sum_{n=0}^{2^N-1} \indicator{\roundaway{SRFF}(x,n)}
      - x  \\
  &=&
   2^{-N}\sum_{n=0}^{2^N-1}
    \indicator{x + 2^{-N} n \ge 1}
    - x
\end{eqnarray*}
and its expectation over the $2^D$ input values is
\begin{eqnarray*}
  bias_{SRFF,D} \kern -25mm &&\\
    &=& 2^{-D} \sum_{i=0}^{2^D-1} bias_{SRFF}(x := 2^{-D} i)  \\
    &=&
    2^{-D} \sum_{i=0}^{2^D-1}  \left( 2^{-N}\sum_{n=0}^{2^N-1}
 \indicator{2^{-D} i + 2^{-N} n \ge 1} - 2^{-D} i
    \right)\\
    &=&
    2^{-D} \sum_{i=0}^{2^D-1}  \left( 2^{-N}\sum_{n=0}^{2^N-1}
 \indicator{2^{-D} i + 2^{-N} n \ge 1} - 2^{-D} i
    \right)\\
    &=&
    \begin{multlined}
      \left(     2^{-D} \sum_{i=0}^{2^D-1}  2^{-N}\sum_{n=0}^{2^N-1}
 \indicator{2^{-D} i + 2^{-N} n \ge 1}
    \right)\\
     - 2^{-D} \sum_{i=0}^{2^D-1} 2^{-D} i
    \end{multlined}
      \\
    &=&
    \begin{multlined}
    \left(    2^{-N}\sum_{n=0}^{2^N-1} 2^{-D} \sum_{i=0}^{2^D-1}
 \indicator{i \ge 2^D - 2^{D-N} n}
    \right)\\ - \frac{1-2^{-D}}2
    \end{multlined} \\
  &=&
  \begin{multlined}
    \left(
        2^{-N}\sum_{n=0}^{2^N-1}
    2^{-D} \sum_{i=\ceil{2^{D} - 2^{D-N} n}}^{2^D-1} 1
    \right)
     - \frac{1-2^{-D}}2
    \end{multlined} \\
    &=&
    \begin{multlined}
      2^{-N}\sum_{n=0}^{2^N-1}
    2^{-D} \biggl(
      (2^D-1)
      - \ceil{2^{D} - 2^{D-N} n}
      + 1\biggr)\\
      - \frac{1-2^{-D}}2    
    \end{multlined} \\
    &=&
    2^{-N}\sum_{n=0}^{2^N-1}
    2^{-D} \floor{2^{D-N} n}
    - \frac{1-2^{-D}}2
\end{eqnarray*}
If $N \le D$ then $2^{D-N}$ is an integer so the floor operation is a no-op.
If $D < N$ then $2^{D-N}$ is a reciprocal power of two so the argument will be integral for some values of~$n$.
We may obtain a bound using $\floor{x}\le x$, which will be tight for $N\le D$.
\begin{eqnarray*}
bias_{SRFF,D}
&=&
    2^{-N}\sum_{n=0}^{2^N-1} 2^{-D} \floor{2^{D-N} n} - \frac{1-2^{-D}}2 \\
&\le&
    2^{-N}\sum_{n=0}^{2^N-1}  2^{-N} n - \frac{1-2^{-D}}2    \\
&=&
    \frac{1-2^{-N}}2     - \frac{1-2^{-D}}2    \\
&=&
    -2^{-(N+1)} + 2^{-(D+1)}\\
&=&
     \frac{2^{-D} - 2^{-N}}2
\end{eqnarray*}

\subsection{Bias of SRF, finite-precision inputs}
\label{asec:fbbias-srf}
The bias at a single $x$ value is
\begin{eqnarray*}
  bias_{SRF,D}(x) &=&
  2^{-N}\sum_{n=0}^{2^N-1} \indicator{\roundaway{SRF}(x,n)}
      - x  \\
  &=&
   2^{-N}\sum_{n=0}^{2^N-1}
    \indicator{x + \left(n+\frac12\right) \times 2^{-N} \ge 1}
    - x
\end{eqnarray*}
and its expectation over the $2^D$ input values is
{\footnotesize
\begin{eqnarray*}
    && 2^{-D} \sum_{i=0}^{2^D-1} bias_{SRFF,D}(x := 2^{-D} i)  \\
    &=&
    2^{-D} \sum_{i=0}^{2^D-1}  \left( 2^{-N}\sum_{n=0}^{2^N-1}
 \indicator{2^{-D} i + \left(n+\frac12\right) \times 2^{-N} \ge 1} - 2^{-D} i
    \right)\\
    &=&
      2^{-D} \sum_{i=0}^{2^D-1}  \left( 2^{-N}\sum_{n=0}^{2^N-1}
  \indicator{2^{-D} i + \left(n+\frac12\right) \times 2^{-N} \ge 1}
  - 2^{-D} i
    \right)\\
    &=&
    \begin{multlined}
      \left(  2^{-D} \sum_{i=0}^{2^D-1}  2^{-N}\sum_{n=0}^{2^N-1}
  \indicator{2^{-D} i + \left(n+\frac12\right) \times 2^{-N} \ge 1}
    \right)\\
     - 2^{-D} \sum_{i=0}^{2^D-1} 2^{-D} i 
    \end{multlined}\\
    &=&
    \begin{multlined}
      \left(    2^{-N}\sum_{n=0}^{2^N-1} 2^{-D} \sum_{i=0}^{2^D-1}
 \indicator{i \ge 2^D - \left(n+\frac12\right) \times 2^{D-N}}
\right)\\
 - \frac{1-2^{-D}}2 
    \end{multlined}
  \\
  &=&
  \begin{multlined}
    \left(    2^{-N}\sum_{n=0}^{2^N-1}
    2^{-D} \sum_{i=\ceil{2^{D} - \left(n+\frac12\right) \times 2^{D-N}}}^{2^D-1} 1
  \right)
    - \frac{1-2^{-D}}2
  \end{multlined}\\
      &=&
      \begin{multlined}
        2^{-N}\sum_{n=0}^{2^N-1}
    2^{-D} \biggl(
      (2^D-1)
      - \ceil{2^{D} - \left(n+\frac12\right) \times 2^{D-N}}
      + 1\biggr)\\
      - \frac{1-2^{-D}}2
    \end{multlined}\\
    &=&
    2^{-N}\sum_{n=0}^{2^N-1}
    2^{-D} \floor{\left(n+\frac12\right) \times 2^{D-N}}
    - \frac{1-2^{-D}}2
      \\
    &=&
      2^{-N}\sum_{n=0}^{2^N-1}
      2^{-D} \floor{n\times 2^{D-N}+\frac12\times 2^{D-N} }
      - \frac{1-2^{-D}}2
    \\
    &=&
      2^{-N}\sum_{n=0}^{2^N-1}
      2^{-D} \floor{n\times 2^{D-N}+2^{D-N-1} }
      - \frac{1-2^{-D}}2
\end{eqnarray*}
}
As before, the floor operation is a complication.
If $D > N$ then $2^{D-N} n$ and $2^{D-N-1}$ are integers, so the bound using $\floor{x}\le x$  will be tight for $N < D$.

Hence $bias_{SRF,D}=$
\begin{eqnarray*}
&&
2^{-N}\sum_{n=0}^{2^N-1}
2^{-D} \floor{n\times 2^{D-N}+2^{D-N-1}}
- \frac{1-2^{-D}}2 \\
&\le&
    2^{-N}\sum_{n=0}^{2^N-1} \left(n\times 2^{-N}+2^{-N-1}\right) - \frac{1-2^{-D}}2
    \\
    &=&
    2^{-2N}\sum_{n=0}^{2^N-1} n    +2^{-N-1} - \frac{1-2^{-D}}2    \\
   &=&
    2^{-2N}\frac{(2^N-1)2^N}2   +2^{-N-1}  - \frac{1-2^{-D}}2
    \\
     &=&
    \frac{1-2^{-N}}2   +2^{-N-1}  - \frac{1-2^{-D}}2
    \\
     &=&
    -2^{-(N+1)} +2^{-N-1} +  2^{-(D+1)}\\
    &=& 2^{-(D+1)}
\end{eqnarray*}
so the bias depends on $D$ but is independent of $N$ (noting that $N \le D$ is required for the computation to be exact.)

\vfill\pagebreak

\bibliographystyle{IEEEtran}
\bibliography{lit}

\end{document}